\documentclass[12pt,a4paper]{article}
\usepackage{lineno,hyperref}
\usepackage[latin5]{inputenc}
\usepackage{amsmath}
\usepackage{graphicx}
\usepackage{subfigure}
\usepackage{float}
\usepackage{mathtools}
\usepackage{color}
\usepackage{cite}
\usepackage{mathrsfs}
\usepackage[]{geometry}
\usepackage[singlelinecheck=false]{caption} 
\usepackage{graphicx}
\usepackage{cite}
\usepackage{authblk}
\usepackage{float}
\usepackage{color}
\usepackage{graphicx}
\usepackage{epstopdf}
\usepackage{amscd,amsmath,amssymb,amsthm,amsfonts,epsfig,graphics}
\DeclareMathOperator{\sech}{sech}
\title{A collocation method based on extended cubic B-splines for numerical solutions of the Klein-Gordon equation}
\author{Alper Korkmaz$^{a,}$\thanks{akorkmaz@karatekin.edu.tr}, Ozlem Ersoy$^{b}$, Idiris Dag$^{c}$  \\
$^{a}${\scriptsize Department of Mathematics, Çankırı Karatekin University, 18200, Çankırı, Turkey.}\\
$^{b}${\scriptsize Department of Mathematics \& Computer, Eskisehir Osmangazi University, 26480, Eskisehir, Turkey.} \\
$^{c}${\scriptsize Department of Computer Engineering, Eskisehir Osmangazi University, 26480, Eskisehir, Turkey.}}
\begin{document}
\maketitle
\begin{abstract} 
A generalization of classical cubic B-spline functions with a parameter is used as basis in the collocation method. Some initial boundary value problems constructed on the nonlinear Klein-gordon equation are solved by the proposed method for extension various parameters. The coupled system derived as a result of the reduction of the time order of the equation is integrated in time by the Crank-Nicolson method. After linearizing the nonlinear term, the collocation procedure is implemented. Adapting the initial conditions provides a linear iteration system for the fully integration of the equation. The validity of the method is investigated by measuring the maximum errors between analytical and the numerical solutions. The absolute relative changes of the conservation laws describing the energy and the momentum are computed for both problems. 
\end{abstract}

\section{Introduction}
\noindent
In the study, we derive numerical solutions for some initial boundary value problems constructed with the nonlinear Klein-Gordon (NKG) equation of the form
\begin{equation}
u_{tt}-u_{xx}-\varepsilon_1u-\varepsilon_2u^3=0
\end{equation}
where $u=u(x,t)$, and $\varepsilon_1$, $\varepsilon_2$ are real parameters\cite{whitham1}. The main part of the equation containing the derivative terms are just the one dimensional wave operator. The remaining part is the derivative of some potential function. The equation was suggested by Klein and Gordon as a relativistic model for a charged particle in an electromagnetic field\cite{deb1}. The laser pulses in two state media, the torsional waves propagating down a stretched wire in a pendula system, the dislocation in crystals, Josephson junction transmission lines, the propagation in ferromagnetic materials of waves carrying rotations of the direction of magnetization are another implementatiton fields of some particular forms of the NKG. The analysis of the rotating black holes can also be possible by this equation\cite{detweiler1}. The geometrical derivation of the NKG requires some particular gauges and coordinate transformations\cite{galehouse1}. The usual procedure for some first order equations in some particular Hilbert spaces supported with some particular norms can lead the scattering field theory for the NKG\cite{weder1}. Some significant properties such as invariance principle, existence and completenes of the wave operators, the intertwining relations were also proven in that study. There exists a dep relation between the NKG and the Schrödinger equations that both can be converted each other\cite{ablowitz1}. The static external field case was also studied for the Klein-Gordon equation\cite{lundberg1}. Employing eigenfunction expansion yielded some important results as strong as like in the Schrödinger Equation in spectral and scattering theory.

\noindent
However, the existence of the envelope type solitons depends upon the sign of the cubic nonlinear terms correlatively the stability of the KGE's dependency\cite{sharma1}. Dusuel et al.\cite{dusuel1} obtained the conditions for the existence of the compactonlike kink solutions of the NKG. They concluded that the static compacton is stable as the dynamic one is not by observing the numerical simulations.

\noindent
The solitary wave solutions at various forms of NKG are determined by Kim and Hong \cite{kim1} using the auxilary equation method based on the solutions of a particular nonlinear ordinary differential equation. They also give the existence conditions of those solutions covering the relations among the parameters and coefficients in the equation. Solitary wave solutions in kink or bell shapes can be constructed by the extended form of the first kind elliptic sub equation method by manipulating the solutions of a first order ordinary differential equation with sixth degree nonlinear term\cite{huang1}. The method is capable to give also explicit forms of the singular and triangular periodic wave solutions. The modified simple equation method is another efficient method to find the solitary wave solutions from the exact traveling solutions under the condition that the parameters in the equation are with their special values\cite{akter1}.

\noindent
Various types of wave solutions like positive or negative frequency plane waves can also be derived from the solutions of convenient field equations\cite{burt1}. Burt and Reid \cite{burt2} set up the exact formal solution of the nonlinear Klein-Gordon equation from the solutions of the linear one. The Klein-gordon equation has also bound state solutions for different attractive potential types\cite{fleischer1}. The soliton interaction is examined in different perspectives by using numerical algorithms\cite{int1,int2,int3}.

\noindent 
Numerical algorithms are also developed for the numerical solutions of the NKG. The classical finite diference method with the central second difference approximation is used to prove the existence of the bounded solutions of the NKG as $t \rightarrow \infty$ \cite{strauss1}. It is also concluded that the degree of the power term causes to change the numer of the oscillations and the amplitude in the solutions.

\noindent
Jiménez and Vázquez \cite{jimenez1} implement four different explicit finite difference schemes and conclude that the scheme which conserves energy is the most suitable one to integrate the NKG to study the long time behaviours of the solutions. Dehghan's study\cite{dehghan1} emphasises that the collocation method based on thin plate spline-radial basis functions can give sufficiently accurate results while solving the inhomogenous NKG with different degreed nonlinear terms. The Fourier collocation method is also implemented to solve some periodic problems\cite{cao1}. That study focuses also the convergence and stability properties of the proposed method. The numerical solutions of the NKG can also be obtained by the multiquadric quasi interpolation method\cite{sarboland1}.

\noindent
The classical polynomial cubic B-spline collocation and unconditionally stable collocation method are derived for the solutions of some initial boundary value problems for the NKG\cite{rash1,zahra1}. In this study we propose a new collocation algorithm based on the extended definiton of the classical polynomial B-splines, namely extended cubic B-splines to solve some initial boundary value problems for the NKG equation. The nature of these B-splines has some siginificant differences from the other B-splines like calassical polynomial \cite{alp3,alp4}, or exponential B-splines\cite{oz1,oz2}. We observe the effects of change of the extention parameter to the accuracy of the solutions.

\noindent 
The order of the NKG in time can be reduced to one to give a nonlinear coupled system 
\begin{equation}
\begin{array}{l}
v_{t}=u_{xx}+\varepsilon_1u+\varepsilon_2u^3 \\ 
u_{t}=v
\end{array}
\label{SYS}
\end{equation}
by assuming $v=u_t$. The initial data 
\begin{equation*}
\begin{aligned}
u\left( x,0\right) &=f\left( x\right) ,\,\, a\leq x\leq b \\ 
\end{aligned}
\end{equation*}
and homogenous Neumann boundary conditions 
\begin{equation}
\begin{aligned}
u_{x}(a,t)&=u_{x}(b,t)=0, \,\, t>0 \\
v_{x}(a,t)&=v_{x}(b,t)=0, \,\, t>0
\end{aligned}
\end{equation}
are chosen in the finite problem interval $[a,b]$ for the convenience.

\section{Numerical Integration of the NKG equation}
The Crank-Nicolson and the suitable forward finite difference for the time integration of the coupled system (\ref{SYS}) yields
\begin{equation}
\begin{aligned}
\dfrac{v^{n+1}-v^{n}}{\Delta t}&=\dfrac{u_{xx}^{n+1}+u_{xx}^{n}}{2}%
+\varepsilon _{1}\dfrac{u^{n+1}+u^{n}}{2}+\varepsilon _{2}\dfrac{%
(u^{3})^{n+1}+(u^{3})^{n}}{2} \\ 
\dfrac{u^{n+1}-u^{n}}{\Delta t}&=\dfrac{v^{n+1}+v^{n}}{2}%
\end{aligned}
\label{e5}
\end{equation}
where $u^{n+1}=u(x,(n+1)\Delta t)$ and $v^{n+1}=v(x,(n+1)\Delta t)$ represent the solutions of the system at the $(n+1)$th. time level. One should note that $t^{n+1}$ equals $t^{n}+\Delta t$, and $\Delta t$ is the time step length, superscripts $n$ and $n+1$ denote the time levels.

Linearization of the term $(u^{3})^{n+1}$ in (\ref{e5})\ as
\begin{equation*}
(u^{3})^{n+1}=3u^{n+1}(u^{2})^{n}-2(u^{3})^{n}  \label{e6}
\end{equation*}%
gives the time-integrated system as
\begin{equation}
\begin{aligned}
\dfrac{v^{n+1}-v^{n}}{\Delta t}&=\dfrac{u_{xx}^{n+1}+u_{xx}^{n}}{2}%
+\varepsilon _{1}\dfrac{u^{n+1}+u^{n}}{2}+\varepsilon _{2}\dfrac{%
3u^{n+1}(u^{2})^{n}-(u^{3})^{n}}{2} \\ 
\dfrac{u^{n+1}-u^{n}}{\Delta t}&=\dfrac{v^{n+1}+v^{n}}{2}%
\end{aligned}
\label{e7}
\end{equation}%

\begin{figure}[ht]
	\centering
		\includegraphics[scale=0.5]{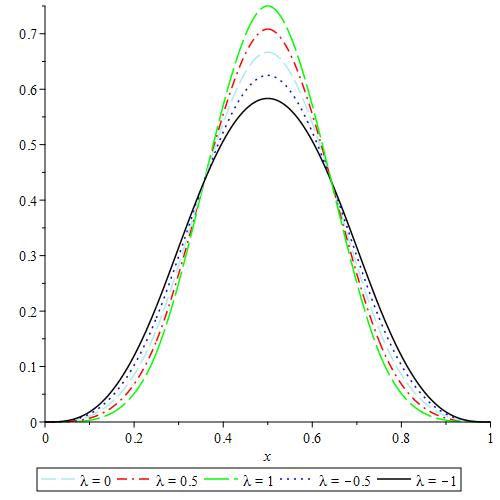}
	\caption{Extended B-splines for various values of the extension parameter $\lambda$}
	\label{fig:Fig1}
\end{figure}

Assume that $\tilde{H}$ is the partition of the finite interval $[a,b]$ as $\tilde{H}:a=x_{1}<x_{2}<\ldots <x_{N}=b$ with equal finite intervals $h=x_i-x_{i-1}, i=2,3,...,N$. An extended cubic B-spline $H_{i}$ is defined as \cite{prenter,dursunext}
\begin{equation}
H_{i}(x)=\frac{1}{24h^{4}}\left \{ 
\begin{array}{ll}
4h(1-\lambda )(x-x_{i-2})^{3}+3\lambda (x-x_{i-2})^{4}, & \left[
x_{i-2},x_{i-1}\right] , \\ 
\begin{array}{l}
(4-\lambda )h^{4}+12h^{3}(x-x_{i-1})+6h^{2}(2+\lambda )(x-x_{i-1})^{2} \\ 
-12h(x-x_{i-1})^{3}-3\lambda (x-x_{i-1})^{4}%
\end{array}
& \left[ x_{i-1},x_{i}\right] , \\ 
\begin{array}{l}
(4-\lambda )h^{4}-12h^{3}(x-x_{i+1})+6h^{2}(2+\lambda )(x-x_{i+1})^{2} \\ 
+12h(x-x_{i+1})^{3}-3\lambda (x-x_{i+1})^{4}%
\end{array}
& \left[ x_{i},x_{i+1}\right] , \\ 
4h(\lambda -1)(x-x_{i+2})^{3}+3\lambda (x-x_{i+2})^{4}, & \left[
x_{i+1},x_{i+2}\right] , \\ 
0 & \text{otherwise.}%
\end{array}%
\right.  \label{e1}
\end{equation}
where $\lambda$ is the real extension parameter. The classical cubic B-spline functions\cite{alp1,alp2} are the particular form when the extension parameter $\lambda$ is chosen as $0$. The set $\{H_i(x) \}_{i=-1}^{N+1}$ defines a basis for the real valued functions defined in the interval $[a,b]$\cite{prenter,dursunext}. The nonzero values of the extension parameter $\lambda$ affects the shape of the cubic B-spline directly, Fig \ref{fig:Fig1}. The relations between the grids and nonzero functional and derivative values of each extended cubic B-spline $H_{i}(x)$ at the partition of $\tilde{H}$ are calculated as in Table \ref{table1}.  

\begin{table}[ht]
	\centering
	\caption{Derivative and functional values of $H_i(x)$ at the grids}
		\begin{tabular}{llllcl}
			\hline \hline
$x$ & $x_{i-2}$ & $x_{i-1}$ & $x_{i}$ & $x_{i+1}$ & \multicolumn{1}{c}{$%
x_{i+2}$} \\ \hline
$24H_{i}(x)$ & $0$ & $4-\lambda $ & $16+2\lambda $ & $4-\lambda $ & 
\multicolumn{1}{c}{$0$} \\ 
$2hH_{i}^{^{\prime }}(x)$ & $0$ & $-1$ & $0$ & $1$ & \multicolumn{1}{c}{$0$} \\ 
$2h^{2}H_{i}^{^{\prime \prime }}(x)$ & $0$ & $2+\lambda $ & $-4-2\lambda $ & $%
2+\lambda $ & \multicolumn{1}{c}{$0$} \\ \hline \hline
		\end{tabular}
	
	\label{table1}
\end{table}

Let $U(x,t)$ and $V(x,t)$ be the approximate solutions to $u(x,t)$ and $v(x,t)$, respectively, defined as
\begin{equation}
\begin{aligned}
U(x,t)&=\sum_{i=-1}^{N+1}\delta _{i}H_{i}(x), \\
V(x,t)&=\sum_{i=-1}^{N+1}\phi _{i}H_{i}(x)  \label{KG1}
\end{aligned}
\end{equation}%
in which $\delta _{i}$ and $\phi _{i}$ are the time dependent new variables. The functional values of $U$ and $V$ and their derivatives are determined by using (\ref{KG1}) as
\begin{equation}
\begin{tabular}{l}
$W_{i}=\dfrac{4-\lambda}{24}\eta _{i-1}+d\frac{8+\lambda}{12}\eta _{i}+\dfrac{4-\lambda}{24}\eta
_{i+1}$ \\ 
$W_{i}^{\prime }=-\dfrac{1}{12h}(\eta _{i-1}-\eta _{i+1})$ \\ 
$W_{i}^{\prime \prime }=\dfrac{2+\lambda}{2h^2}(\eta _{i-1}-2\eta
_{i}+\eta _{i+1})$%
\end{tabular}%
\end{equation}%
where $W_i$ and $\eta$ are general notations. $\eta$ denotes $\delta$ when $W_i$ is chosen as $U_i=U(x_i)$ and it stands for $\phi$ when $W_i=V_i=V(x_i)$.

Substituting the approximate solutions $U_i$ and $V_i$ defined in (\ref{KG1}) and their derivatives into (\ref%
{e7}) and rearranging the resulting equations yields the iterative system
\begin{eqnarray}
&&\omega _{m1}\delta _{m-1}^{n+1}+\omega _{m2}\phi _{m-1}^{n+1}+\omega _{m3}\delta
_{m}^{n+1}+\omega _{m4}\phi _{m}^{n+1}+\omega _{m1}\delta _{m+1}^{n+1}+\omega
_{m2}\phi _{m+1}^{n+1}  \label{e9} \\
&=&\omega _{m5}\delta _{m-1}^{n}+\omega _{m2}\phi _{m-1}^{n}+\omega _{m6}\delta
_{m}^{n}+\omega _{m4}\phi _{m}^{n}+\omega _{m5}\delta _{m+1}^{n}+\omega _{m2}\phi
_{m+1}^{n}  \notag
\end{eqnarray}%
\begin{eqnarray}
&&\omega _{m2}\delta _{m-1}^{n+1}+\omega _{m7}\phi _{m-1}^{n+1}+\omega _{m4}\delta
_{m}^{n+1}+\omega _{m8}\phi _{m}^{n+1}+\omega _{m2}\delta _{m+1}^{n+1}+\omega
_{m7}\phi _{m+1}^{n+1}  \label{e10} \\
&=&\omega _{m2}\delta _{m-1}^{n}-\omega _{m7}\phi _{m-1}^{n}+\omega _{m4}\delta
_{m}^{n}-\omega _{m8}\phi _{m}^{n}+\omega _{m2}\delta _{m+1}^{n}-\omega _{m7}\phi
_{m+1}^{n}  \notag
\end{eqnarray}

The coefficients of equation system (\ref{e9}) and (\ref{e10}) for KGE
equation can be determined as follow%
\begin{equation*}
\begin{array}{l}
\omega _{m1}=\left( -3\varepsilon _{2}K^{2}-\varepsilon _{1}\right) \alpha
_{1}-\gamma _{1} \\ 
\omega _{m2}=\dfrac{2}{\Delta t}\alpha _{1} \\ 
\omega _{m3}=\left( -3\varepsilon _{2}K^{2}-\varepsilon _{1}\right) \alpha
_{2}-\gamma _{2} \\ 
\omega _{m4}=\dfrac{2}{\Delta t}\alpha _{2} \\ 
\omega _{m5}=\left( \varepsilon _{1}-\varepsilon _{2}K^{2}\right) \alpha
_{1}+\gamma _{1} \\ 
\omega _{m6}=\left( \varepsilon _{1}-\varepsilon _{2}K^{2}\right) \alpha
_{2}+\gamma _{2} \\ 
\omega _{m7}=-\alpha _{1} \\ 
\omega _{m8}=-\alpha _{2}%
\end{array}%
\end{equation*}

where
\begin{equation*}
\begin{aligned}
K&=\alpha _{1}\delta _{i-1}^n+\alpha _{2}\delta _{i}^n+\alpha _{1}\delta _{i+1}^n \\
\alpha_1&=\frac{4-\lambda}{24} \\
\alpha_2&=\frac{8+\lambda}{12} \\
\gamma_1&=\frac{2+\lambda}{2h^2} \\
\gamma_2&=-\frac{4+2\lambda}{2h^2}
\end{aligned}
\end{equation*}

\noindent
The system (\ref{e9}) and (\ref{e10}) can be rewritten in the matrix notation for the sake of simplicity as
\begin{equation}
\mathbf{Ax}^{n+1}=\mathbf{Bx}^{n}  \label{e11}
\end{equation}%
where%
\begin{equation*}
\mathbf{A=}%
\begin{bmatrix}
\omega _{m1} & \omega _{m2} & \omega _{m3} & \omega _{m4} & \omega _{m1} & \omega _{m2} &  & 
&  &  \\ 
\omega _{m2} & \omega _{m7} & \omega _{m4} & \omega _{m8} & \omega _{m2} & \omega _{m7} &  & 
&  &  \\ 
&  & \omega _{m1} & \omega _{m2} & \omega _{m3} & \omega _{m4} & \omega _{m1} & \omega _{m2} & 
&  \\ 
&  & \omega _{m2} & \omega _{m7} & \omega _{m4} & \omega _{m8} & \omega _{m2} & \omega _{m7} & 
&  \\ 
&  &  & \ddots & \ddots & \ddots & \ddots & \ddots & \ddots &  \\ 
&  &  &  & \omega _{m1} & \omega _{m2} & \omega _{m3} & \omega _{m4} & \omega _{m1} & \omega
_{m2} \\ 
&  &  &  & \omega _{m2} & \omega _{m7} & \omega _{m4} & \omega _{m8} & \omega _{m2} & \omega
_{m7}%
\end{bmatrix}%
\end{equation*}%
and%
\begin{equation*}
\mathbf{B=}%
\begin{bmatrix}
\omega _{m5} & \omega _{m2} & \omega _{m6} & \omega _{m4} & \omega _{m5} & \omega _{m2} &  & 
&  &  \\ 
\omega _{m2} & -\omega _{m7} & \omega _{m4} & -\omega _{m8} & \omega _{m2} & -\omega _{m7} & 
&  &  &  \\ 
&  & \omega _{m5} & \omega _{m2} & \omega _{m6} & \omega _{m4} & \omega _{m5} & \omega _{m2} & 
&  \\ 
&  & \omega _{m2} & -\omega _{m7} & \omega _{m4} & -\omega _{m8} & \omega _{m2} & -\omega _{m7}
&  &  \\ 
&  &  & \ddots & \ddots & \ddots & \ddots & \ddots & \ddots &  \\ 
&  &  &  & \omega _{m5} & \omega _{m2} & \omega _{m6} & \omega _{m4} & \omega _{m5} & \omega
_{m2} \\ 
&  &  &  & \omega _{m2} & -\omega _{m7} & \omega _{m4} & -\omega _{m8} & \omega _{m2} & 
-\omega _{m7}%
\end{bmatrix}%
\end{equation*}

The system (\ref{e11}) has $2N+2$ linear equations and $2N+6$ unknown parameters described as $%
\mathbf{x}^{n+1}=(\delta _{-1}^{n+1},\phi _{-1}^{n+1},\delta _{0}^{n+1},\phi
_{0}^{n+1}\ldots ,\delta _{N+1}^{n+1},\phi _{N+1}^{n+1})$. A unique solution
of this system requires the equal number of equations and parameters. Implement of the boundary conditions $$U_{x}(a,t)=0,U_{x}(b,t)=0,V_{x}(a,t)=0,V_{x}(b,t)=0$$ equalize the number of unknown parameters by generating relations 
\begin{equation*}
\delta _{-1}=\delta _{1},\text{ }\phi _{-1}=\phi _{1},\text{ }\delta
_{N-1}=\delta _{N+1},\text{ }\phi _{N-1}=\phi _{N+1}
\end{equation*}

When the parameters $\delta _{-1},\phi _{-1},\delta _{N+1},\phi
_{N+1}$ are eliminated from the system, we have $2N+2$ linear equations with $2N+2$ unknowns. We solve this system of linear equation by the Thomas algorithm for the systems having six-banded coefficient matrices that is adapted from the algorithm for the systems having seven-banded coefficient matrix. In order to start the iteration algorithm, we need the initial vector $\mathbf{x^0}$.  Assuming $\mathbf{x_{1}^0}=(\delta _{-1}^0,\delta _{0}^0,..\delta _{N}^0,\delta _{N+1}^0)$, $\mathbf{x_{2}^0}=(\phi
_{-1}^0,\phi _{0}^0,..\phi _{N}^0,\phi _{N+1}^0)$ are the components of the initial vector $\mathbf{x^0}$ of the iteration, the parameters are eliminated by using the equalities 
\begin{equation}
\begin{aligned}
U_{x}(a,0)&=0=\delta _{-1}^{0}-\delta _{1}^{0}, \\ 
U_{x}(x_{i},0)&=\delta _{i-1}^{0}-\delta _{i+1}^{0}=U_{x}(x_{i},0),i=1,...,N-1
\\ 
U_{x}(b,0)&=0=\delta _{N-1}^{0}-\delta _{N+1}^{0}, \\ 
V_{x}(a,0)&=0=\phi _{-1}^{0}-\phi _{1}^{0} \\ 
V_{x}(x_{i},0)&=\phi _{i-1}^{0}-\phi _{i+1}^{0}=V_{x}(x_{i},0),i=1,...,N-1 \\ 
V_{x}(b,0)&=0=\phi _{N-1}^{0}-\phi _{N+1}^{0}%
\end{aligned}%
\end{equation}
to be able to start the iteration (\ref{e11}).
\section{Numerical Solutions}

\noindent
This section is devoted to focus the perform of the suggested method by implementing it to some initial boundary value problems for the NKG. The discrete maximum error norm
\begin{equation*}
\begin{aligned}
L_{\infty}(t)&=\left \vert u(x,t)-U(x,t)\right \vert _{\infty }=\max \limits_{i}\left
\vert u(x_i,t)-U(x_i,t)\right \vert \\
\end{aligned}%
\end{equation*}
is defined to check the validity and accuracy of the suggested method by measuring the error between the analytical and numerical solution at a specific time $t$. The conservation of the energy(E) and the momentum(P) defined as \cite{whitham1,deb1,jhangeer1}
\begin{equation}
\begin{aligned}
E&=\frac{1}{2}\int\limits_{-\infty}^{\infty}{u_t^2+u_x^2-\varepsilon_1u^2-\frac{1}{2}\varepsilon_2u^4dx} \\
P&=\int\limits_{-\infty}^{\infty}{u_xu_tdx}
\end{aligned}
\end{equation} 
can also be a good indicator of an efficient method in case the absence of the analytical solutions. We define absolute relative changes $C(E_t)$ and $C(P_t)$ at the time $t$ of the conserved quantities $E$ and $P$ as
\begin{equation}
\begin{aligned}
C(E_t)&=\left | \frac{E_t-E_0}{E_0} \right | \\
C(P_t)&=\left| \frac{P_t-P_0}{P_0} \right|
\end{aligned}
\end{equation}
where $E_0$ and $P_0$ are initial values of the energy and the momentum of the system, respectively.

\subsection{Traveling Wave Case}

\noindent The initial boundary value problem is defined for $\varepsilon_1=1$, $\varepsilon_2=-1$ in the NKG equation. The analytical solution 
\begin{equation}
u(x,t)=\tanh (\frac{(x-\nu t)}{\sqrt{2(1-\nu^{2})}})  \label{n1}
\end{equation}%
describes a traveling wave moving along the $x-$axis with the velocity $|\nu|<1$\cite{zaki1}. The initial data are derived from the analytical solution (\ref{n1}) by substituting $t=0$ into it. The Neumann conditions at both ends of the interval $[-30,30]$ are used for the numerical solutions. The routine is run for various values of the discretization parameters $h$ and $\Delta t$ with the choice of the velocity as $\nu=0.5$ to the time $t=10$. In order to improve the accuracy of the results, the extension parameter is scanned between $[-1,1]$ with the increment $\Delta \lambda = 0.0001$ for the optimum choice of the extension parameter.

\begin{figure}[h]
    \subfigure[Traveling wave simulation]{
   \includegraphics[scale =0.6] {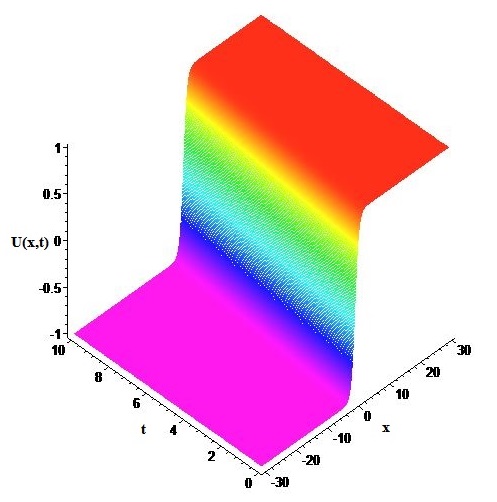}
   \label{fig:1a}
 }
 \subfigure[The maximum error distribution]{
   \includegraphics[scale =0.6] {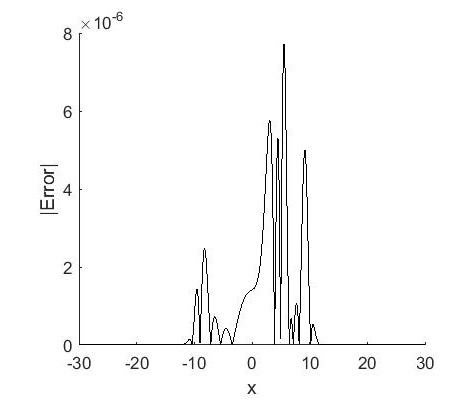}
   \label{fig:1b}
 }
 \caption{Traveling wave simulation and the maximum error distribution at $t=10$}
\end{figure}

\noindent
The initial data is an $S-$shaped wave positioned at the origin with the properties $\lim\limits_{x\rightarrow -\infty} {U(x,t)}\rightarrow -1$, $\lim\limits_{x\rightarrow \infty} {U(x,t)}\rightarrow 1 $. When the simulation starts, the wave moves back along the $x-$axis with the constant velocity $\nu=0.5$ without changing its shape, Fig \ref{fig:1a}. 

\noindent
The maximum error distribution for the optimum value of the extension parameter $\lambda$ and the discretization paramters $h=0.02$ and $\Delta t=0.005$ at the time $t=10$ is depicted in Fig (\ref{fig:1b}). It is observed from both figures that the error accumulates at the points where the wave descent occours.

\noindent
The discrete maximum norms for both $\lambda =0$ and various $\lambda$ cases are tabulated in Table \ref{tab:3}. When the discretization parameters are chosen as $h=0.2$ and $\Delta t=0.05$, the maximum error is in two decimal digits for $\lambda=0$. This error is improved to four decimal digit accuracy by choosing the optimum extension parameter as $\lambda=-0.0101$. Reducing the dicsretization parameters to $h=0.1$ and $\Delta t=0.02$ improves the results five times for both $\lambda =0$ and optimum $\lambda =-0.0026$. The discretization parameters $h=0.05$ and $\Delta t=0.01$ gives four decimal digit accurate results for $\lambda =0$. The maximum error is reduced to the fifth decimal in this case by determining the optimum extension parameter as $\lambda =-0.007$. One more reduce of the discretization parameters to $h=0.02$ and $\Delta t=0.005$ gives the error in the fourth decimal digit for $\lambda=0$ as provides six decimal digits accuracy in the results for $\lambda =-0.0001$. 

\begin{table}[H]
\caption{Discrete maximum norms for various discretization parameters $h$ and $\Delta t$}  
\begin{tabular}{llll}
\hline \hline
$h$ & $\Delta t$ & $L_{\infty}(10)(\lambda=0)$ &  $L_{\infty}(10)$ \\
\hline
$0.2$ & $0.05$  & $1.0709\times 10^{-2}$ & $6.4162\times 10^{-4}$ ($\lambda = -0.0101$) \\ 
$0.1$ & $0.02$  &  $2.7968\times 10^{-3}$ & $1.6618\times 10^{-4}$ ($\lambda = -0.0026$) \\ 
$0.05$ & $0.01$ &  $7.0161\times 10^{-4}$ & $8.1010\times 10^{-5}$ ($\lambda = -0.0007$) \\ 
$0.02$ & $0.005$ &  $1.0984\times 10^{-4}$& $7.7240\times 10^{-6}$ ($\lambda = -0.0001$) \\ \hline \hline
\end{tabular}%
\label{tab:3}
\end{table}

\noindent
The initial values of the conservation laws are determined by using symbolic software as
\begin{equation}
\begin{aligned}
E_0&=-\frac{1}{9}\,{\frac {405\,{{\rm e}^{20\,\sqrt {2}\sqrt {3}}}+405\,{{\rm e}^{
40\,\sqrt {2}\sqrt {3}}}+12\,{{\rm e}^{20\,\sqrt {2}\sqrt {3}}}\sqrt {
2}\sqrt {3}+\tilde{A}
}{1+3\,{{\rm e}^{
20\,\sqrt {2}\sqrt {3}}}+3\,{{\rm e}^{40\,\sqrt {2}\sqrt {3}}}+{
{\rm e}^{60\,\sqrt {2}\sqrt {3}}}}} \\
\tilde{A}&=4\,\sqrt {2}\sqrt {3}-12\,{{\rm e}^{40\,\sqrt {2}\sqrt {3}
}}\sqrt {2}\sqrt {3}-4\,{{\rm e}^{60\,\sqrt {2}\sqrt {3}}}\sqrt {2}
\sqrt {3}+135\,{{\rm e}^{60\,\sqrt {2}\sqrt {3}}}+135 \\
P_0&=-\frac{2}{9}\,\sqrt {2}\sqrt {3}
\end{aligned}
\end{equation} 
The approximate values of the conservation laws are computed as $E_0=-13.91133789$ and $P_0=-0.5443310539$ initially. It should be noted that the initial quantity of the energy $E_0$ is computed by reducing the bounds of the related integral to the problem interval $[-30,30]$. The absolute relative changes of both conserved quantities are reported at the simulation terminating time $t=10$ in Table \ref{tab:4}. The absolute relative change of the energy of the system is in seven decimal digits and of the momentum is in five decimal digits for both $\lambda =0$ and the optimum $\lambda=-0.0101$ when the discretization parameters are $h=0.2$ and $\Delta t=0.05$. The absolute relative changes of the energy and the momentum are in eight and six decimal digits, respectively for $h=0.1$ and $\Delta t=0.02$. Reducing the discretization parameters to $h=0.05$ and $\Delta t=0.01$ gives nine decimal digits absolute relative change for the energy and seven decimal digits for the momentum. When $h=0.02$ and $\Delta t=0.005$, the absolute relative changes are in ten and eight decimal digits for the energy and the momentum, respectively. The reduction of the discretization parameters improves the absolute relative changes of both the energy and momentum quantities but we do not observe a significant improve on the results with respect to the optimum choice of the extension parameter $\lambda$.

\begin{table}[H]
\caption{Discrete maximum norms for various discretization parameters $h$ and $\Delta t$}  
\begin{tabular}{llllll}
\hline \hline
$h$ & $\Delta t$ & $C(E_{10})$ & $C(E_{10})$ & $C(P_{10})$& $C(P_{10})$  \\
&&$\lambda =0$&optimum $\lambda$&$\lambda =0$&optimum $\lambda$\\ 
\hline
$0.2$ & $0.05$ & $2.4670\times 10^{-7}$ & $4.3656\times 10^{-7}$ & $3.8299\times 10^{-5}$& $3.6371\times 10^{-5}$  \\ 
$0.1$ & $0.02$ & $1.7897\times 10^{-8}$& $2.8211\times 10^{-8}$ & $2.4033\times 10^{-6}$& $2.3174\times 10^{-6}$   \\ 
$0.05$ & $0.01$ & $1.6705\times 10^{-9}$& $4.5956\times 10^{-9}$ & $2.9338\times 10^{-7}$& $2.8741\times 10^{-7}$   \\ 
$0.02$ & $0.005$ & $6.0376\times 10^{-10}$& $6.7339\times 10^{-10}$ & $3.6019\times 10^{-8}$ & $3.5804\times 10^{-8}$   \\ 
\hline \hline
\end{tabular}%
\label{tab:4}
\end{table}
\subsection{Single Solitary Wave Case}
\begin{figure}[h]
    \subfigure[Propagation of a single solitary wave]{
   \includegraphics[scale =0.6] {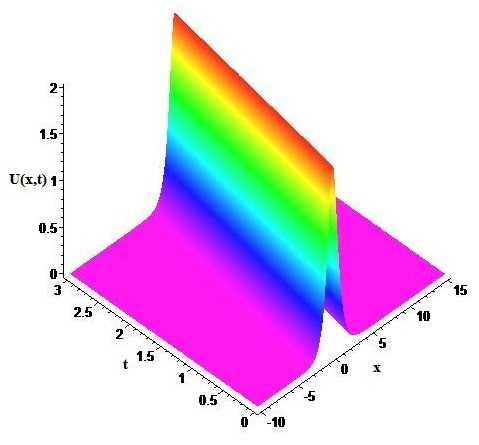}
   \label{fig:2a}
 }
 \subfigure[The maximum error distribution]{
   \includegraphics[scale =0.5] {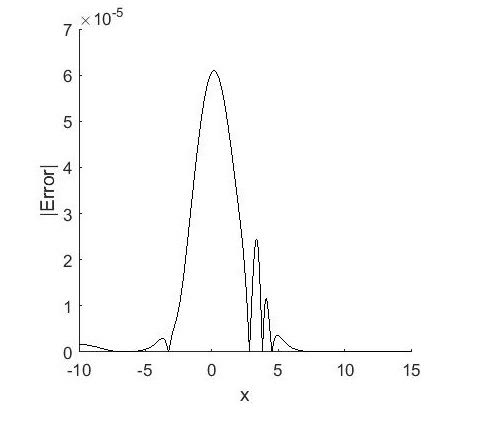}
   \label{fig:2b}
 }
 \caption{Propagation of a single solitary wave and the maximum error distribution at $t=3$}
\end{figure}

The single solitary wave solution of the NKG is derived from the solution in Polyanin's book\cite{polyanin1} as
\begin{equation}
u(x,t)=2\sech{(\sqrt{2}\{\sinh{(1)}\}x-\{ \cosh{(1)}\}t)} \label{ana1}
\end{equation}
for $\varepsilon_1=2$, $\varepsilon_2=-1$. The solution models the propagation of a single solitary wave of amplitude $2$ to the right along the horizontal axis. The peak of the wave is positioned at the origin initially. The problem interval is shrunk to $[-10,15]$ to be able apply the numerical method. The initial data is obtained by substituting $t=0$ in the analytical solution (\ref{ana1}). The suitable Neumann conditions are considered in accordance with the analytical solution. The routine is run up to the time $t=3$ for various values of the discretization parameters. The simulation of the motion and the maximum error distribution are depicted in Fig \ref{fig:2a} and Fig \ref{fig:2b}. The peak is positioned at $x=1.31$ when the simulation time reaches $t=1$. The position of the peak is measured as $x=2.625$ at $t=2$ and as $x=3.94$ at $t=3$. Thus, the average velocity of the wave can be computed approximately as $1.32$ in the appropriate units. 

\noindent
Even though we scan the extension parameter in the interval $[-1,1]$ with the step size $\Delta \lambda =0.0001$, the best results are obtained when the extension parameter is zero. The results are summerized in Table \ref{tab:5}. The discrete maximum error is measured in four decimal digit accuracy for $h=0.05$ and $\Delta t=0.01$ at $t=1$. The maximum errors are in three decimal digits at the times $t=2$ and $t=3$ with the same discretization parameters. When the time step size is reduced to $0.001$, the results are worse than the results obtained with $\Delta t=0.01$ at $t=1$. Even though the accuracy decimals are equal at $t=2$ and at $t=3$, the decimal values of the accuracy of the results are worse than the ones obtained for $\Delta t=0.01$. When $h$ is reduced ten times, an improve is observable in the results. The maximum errors are determined as $3.8323\times 10^{-4}$,  $5.6663\times 10^{-4}$ and $1.5227\times 10^{-3}$ at the times $t=1$, $t=2$ and $t=3$, respectively with the time step size $\Delta t=0.01$. When $h=0.005$ and $\Delta t=0.001$, the accuracy of the results are in six decimal digits at $t=1$, and five decimal digits at $t=2$ and $t=3$.
\begin{table}[H]
\caption{Discrete maximum norms for various discretization parameters and $\lambda=0$}  
\begin{tabular}{llllll}
\hline \hline
$h$ & $\Delta t$ & $L_{\infty }(1)$ & $L_{\infty }(2)$ & $L_{\infty }(3)$ 
\\ 
$0.05$ & $0.01$ & $8.5481\times 10^{-4}$ & $1.9112\times 10^{-3}$ & $%
6.0948\times 10^{-3}$   \\ 
& $0.001$ & $1.2412\times 10^{-3}$ & $2.4715\times 10^{-3}$ & $4.7686\times
10^{-3}$  \\ 
$0.005$ & $0.01$ & $3.8323\times 10^{-4}$ & $5.6663\times 10^{-4}$ & $%
1.5227\times 10^{-3}$  \\ 
& $0.001$ & $8.5593\times 10^{-6}$ & $1.9112\times 10^{-5}$ & $6.0943\times
10^{-5}$   \\
 \hline \hline
\end{tabular}%
\label{tab:5}
\end{table}

\noindent
The conserved quantities describing the energy and momentum of the system is calculated using symbolic calculation software as
\begin{equation}
\begin{aligned}
E_0&=\frac{8}{3}\frac{\sqrt{2}(\cosh^2{1}-1)}{\sinh{1}}\\ 
P_0&=-\frac{8}{3}\sqrt{2}\cosh{1}
\end{aligned}
\end{equation}
with the approximate values $E_0=4.431961243$ and $P_0=-5.819321497$ initially. The absolute relative changes of these two quantities are tabularised in Table \ref{tab:6}. The absolute relatives change in the energy are measured in six decimal digits when $h=0.05$, and $h=0.005$ with $\Delta t=0.01$ at the time $t=3$. Reducing $\Delta t$ to $0.001$ improves the change to nine decimal digits for the absolute relative change of the energy. The absolute relative change of the momentum is in six decimal digits when $\Delta t=0.01$ for both $h=0.05$ and $h=0.005$. Choosing $h=0.05$ and $\Delta t=0.001$ gives eight decimal digit absolute relative change as gives nine decimal digit absolute relative change when $h=0.005$ and $\Delta t=0.001$.
\begin{table}[H]
\caption{Absolute relative changes of the conserved quantitites for $\lambda =0$}  
\begin{tabular}{llll}
\hline \hline
$h$ & $\Delta t$ & $C(E_{3})$ & $C(P_{3})$  \\ 
\hline
$0.05$ & $0.01$ & $4.7265\times 10^{-6}$ & $1.5699\times 10^{-6}$ \\ 
& $0.001$ & $2.7641\times 10^{-6}$ & $1.1944\times 10^{-8}$  \\ 
$0.005$ & $0.01$ & $4.2897\times 10^{-6}$ & $2.5049\times 10^{-6}$  \\ 
& $0.001$ & $4.2706\times 10^{-9}$ & $2.3602\times 10^{-9}$  \\
\hline \hline
\end{tabular}%
\label{tab:6}
\end{table}

\section{Conclusion}
\noindent
The extended form of the cubic polynomial B-splines are used as basis in the collocation method for the solutions of the nonlinear Klein-Gordon equation. The order of the NKG is reduced to one to be able to integrate in time by Crank-Nicolson method. The dependent variables in the resulting system are approximated by the extended cubic B-splines. The validity and accuracy of the suggested method are by solving two initial boundary value problems. The discrete maximum error norms and absolute relative changes of the conserved quantities are reported to validate the results.

\noindent
The first problem describing the travel of a $\tanh-$type wave is solved succesfully by the suggested method. The scan of the extension parameter improves the results when compared with the results of the classical polynomial cubic B-spline case. 

\noindent
In the second problem, we study the propagation of a single solitary wave. The numerical results are in a good agreement with the analytical ones. In contrast to the first example, the scan of the extension parameter does not improve the results in this case.

\noindent
The absolute relative chances of the conserved quantities correspond the theoretical aspects of the conservation laws.

Acknowledgements: \textsl{A brief part of this study was presented orally in International Conference on Applied Mathematics and Analysis, Ankara-Turkey, 2016.}

\end{document}